\documentclass[twoside,11pt,leqno]{article}
\usepackage{amsfonts}

 \def\X{{\cal X}}  \def\H{{\cal H}} 
\def\d{\rm{des}}

\def\P{{\bf {P}}} \def\E{\cal{E}} 
\def\A{\mathbb{A}} \def\e{\mathbb{B}} \def\L{\mathbb{L}} \def\R{\mathbb{R}}
\def\S{\mathbb{S}} \def\T{\mathbb{T}} \def\d{\sum_{i=1}^d} \def\j{\sum_{j=0}^m}

\def\I{\mathbb{I}}
\def\B{B({\cal H})} \def\b{B({\cal X})}

\newtheorem{df}{Definition}[section]
\newtheorem{thm}[df]{Theorem} \newtheorem{pro}[df]{Proposition}
 
\newtheorem{rema}[df] {Remark} \newtheorem{lem}[df] {Lemma}
\def\sfstp{{\hskip-1em}{\bf.}{\hskip1em}}

\def\subject#1{\renewcommand{\thefootnote}{}\footnote
{AMS(MOS) subject classification (2010). Primary: {#1}}}

\def\keywords#1{\renewcommand{\thefootnote}{}\footnote
{Keywords: {#1}}}

\def\enddemo{\qed \endtrivlist} \expandafter\let\csname
enddemo*\endcsname=\enddemo

\def\qedsymbol{\ifmmode\bgroup\else$\bgroup\aftergroup$\fi
\vcenter{\hrule\hbox{\vrule
height.5em\kern.5em\vrule}\hrule}\egroup}
\def\qed{\ifmmode\else\unskip\nobreak\fi\quad\qedsymbol}

\title{\bf On strict isometric and strict symmetric  commuting $d$-tuples of Banach space operators}
\author{\normalsize B.P.~Duggal}
\date{}

\begin{document}

\maketitle \thispagestyle{empty} \vskip-16pt

\subject{47A05, 47A55; Secondary47A11, 47B47.} \keywords{ Banach space,  commuting $d$-tuples,left/right multiplication operator,  $m$-left invertible, $m$-isometric and $m$-selfadjoint operators, product of operators, tensor product .  }

\begin{abstract} Given commuting $d$-tuples $\S_i$ and $\T_i$, $1\leq i\leq 2$, Banach space operators such that the tensor products pair $(\S_1\otimes\S_2,\T_1\otimes\T_2)$ is strict $m$-isometric (resp., $\S_1$,  $\S_2$ are invertible and  $(\S_1\otimes\S_2,\T_1\otimes\T_2)$ is strict $m$-symmetric), there exist integers $m_i >0$, and a non-zero scalar $c$, such that $m=m_1+m_2-1$, $(\S_1,{\frac{1}{c}}\T_1)$ is strict $m_1$-isometric and $(\S_2, c\T_2)$ is strict $m_2$-isometric (resp., there exist integers $m_i >0$, and a non-zero scalar $c$, such that $m=m_1+m_2-1$, $(\S_1,{\frac{1}{c}}\T_1)$ is strict $m_1$-symmetric and $(\S_2, c\T_2)$ is strict $m_2$-symmetric. However, $(\S_i,\T_i)$ is strict $m_i$-isometric (resp., strict $m_i$-symmetric) for $1\leq i\leq 2$ implies only that  $(\S_1\otimes\S_2,\T_1\otimes\T_2)$ is  $m$-isometric (resp.,  $(\S_1\otimes\S_2,\T_1\otimes\T_2)$ is  $m$-symmetric).
\end{abstract}


\section {\sfstp Introduction} Let $\b$ (resp., $\B$) denote the algebra of operators, i.e. bounded linear transformations, on an infinite dimensional complex Banach space $\X$ into itself (resp., on an infinite dimensional complex Hilbert space $\H$ into itself) , $\mathbb C$ denote the complex plane, $\b^d$ (resp., $\B^d$ and  ${\mathbb C}^d$) the product of $d$ copies of $\b$ (resp., $\B$ and $\mathbb C$) for some integer $d\geq 1$, $\overline{z}$ the conjugate of $z\in\mathbb C$ and ${\bf z}=(z_1,z_2,...,z_d)\in{\mathbb C}^d$. A $d$-tuples $\A=(A_1, \cdots, A_d)\in\b^d$ is a commuting $d$-tuple if $[A_i,A_j]=A_iA_j-A_jA_i=0$ for all $1\leq i,j\leq d$. If $\P$ is a polynomial in ${\mathbb C}^d$ and ${\A}$ is a $d$-tuple of commuting operators in $\B^d$, then  ${\A}$ is a hereditary root of $\P$ if $\P({\A})=0$. Two particular operator classes of hereditary roots which have drawn a lot of attention in the recent past are those of  $m$-isometric and $m$-symmetric (also called $m$-selfadjoint) operators, where $A\in\B$ is $m$-isometric, $m$ some positive integer, if $\sum_{j=0}^m(-1)^j\left(\begin{array}{clcr}m\\j\end{array}\right)A^{*j}A^j=0$ and $A\in\B$ is $m$-symmetric if $\sum_{j=0}^m(-1)^j\left(\begin{array}{clcr}m\\j\end{array}\right){A^*}^{(m-j)}A^j=0.$
Clearly, $m$-isometric operators arise as  solutions of $P(z)=({\overline{z}}z-1)^m=0$ and   $m$-symmetric operators  arise as solutions of $P(z)=({\overline{z}}-z)^m=0$. The class of $m$-isometric operators was introduced by Agler \cite{A} and the class of $m$-symmetric operators was introduced by Helton \cite{Hel} ({\it albeit} not as operator solutions of the polynomial equation $({\overline{z}}-z)^m=0$). These classes of operators, and their variants (of the {\em left $m$-invertible and $m$-symmetric Banach space pairs $(A,B)$ type} \cite{DK2}),  have since been studied by a multitude of authors, amongst them Agler and Stankus \cite{AS1}, Bayart \cite{FB}, Bermudez {\it et al} \cite{BMN},  Duggal and Kim \cite{{DK1}, {DK2}},  Gu \cite{G},  Gu and Stankus \cite{GS}  and Paul and Gu \cite{PG}.

\

Generalising the $m$-isometric property of operators $A\in\B$ to commuting $d$-tuples $\A\in \B^d$, Gleeson and Richter \cite {GR} say that $\A$ is $m$-isometric if 
\begin{eqnarray}
\sum_{j=0}^m{(-1)^j\left(\begin{array}{clcr}m\\j\end{array}\right)\sum_{|\beta|=j}{\frac{j!}{{\beta}!}}\A^{*\beta}\A^{\beta}}=0,
\end{eqnarray}
where
$$
\beta=(\beta_1, \cdots, \beta_d),  \ |\beta|=\sum_{i=1}^d \beta_i, \ {\beta}!=\Pi_{i=1}^d{\beta_i}!, \A^{\beta}=\Pi_{i=1}^d{A_i^{\beta_i}}, \A^{*{\beta}}=\Pi_{i=1}^d{A_i^{*{\beta_i}}};
$$
$\A$ is said to be $m$-symmetric if 
\begin{eqnarray}
\sum_{j=0}^m{(-1)^j\left(\begin{array}{clcr}m\\j\end{array}\right)(A_1^*+ \cdots +A^*_d)^{m-j}(A_1+ \cdots +A_d)^j}=0.
\end{eqnarray}
These generalisations and certain of their variants, in particular  left $m$-invertible Banach space operator pairs $(A,B)$: 
$$
\triangle^m_{A,B}(I)=\sum_{j=0}^m{(-1)^j\left(\begin{array}{clcr}m\\j\end{array}\right) A^jB^j}=0
$$ 
and $m$-symmetric pairs $(A,B)$:
$$\delta^m_{A,B}(I)=\sum_{j=0}^m{(-1)^j\left(\begin{array}{clcr}m\\j\end{array}\right)A^{m-j}B^j}=0,$$ 
have recently been the subject matter of a number of studies, see \cite{{CA}, {DK2}, {DM}, {G}, {PG}, {ACL}} for further references. 

\

Recall that a pair $(A,B)$ of Banach space operators is strict $m$-left invertible if $\triangle^m_{A,B}(I)=0$ and $\triangle^{m-1}_{A,B}(I)\neq 0$; similarly, the pair $(A,B)$ is strict $m$-symmetric if $\delta^m_{A,B}(I)=0$ and $\delta^{m-1}_{A,B}(I)\neq 0$. Products $(A_1A_2,B_1B_2)$ of $m_i$-isometric, similarly $m_i$symmetric, pairs $(A_i,B_i)$, $1\leq i\leq 2$, such that $A_1$ commutes with $A_2$, and $B_1$ commutes with $B_2$, are $(m_1+m_2-1)$-isometric, respectively $(m_1+m_2-1)$-symmetric \cite{{BMN}, {DM}, {DK2}, {G}}. The converse fails, even for  strict $m$-isometric (and strict $m$-symmetric) operator pairs $(A_1A_2,B_1B_2)$. A case where there is an answer in the positive is that of the tensor product pairs $(A_1\otimes A_2, B_1\otimes B_2)$: (i) $(A_1\otimes A_2,B_1\otimes B_2)$ is  $m$-isometric if and only if there exist integers $m_i >0$, and a non-zero scalar $c$, such that $m=m_1+m_2-1$, $(A_1,{\frac{1}{c}} B_1)$ is  $m_1$-isometric and $(A_2, c B_2)$ is  $m_2$-isometric \cite[Theorem 1.1]{PG}; (ii) if $A_i$ are left invertible and $B_i$ are right invertible, $1\leq i\leq 2$, then  $(A_1\otimes A_2,B_1\otimes B_2)$ is  $m$-symmetric if and only if there exist integers $m_i >0$ and a non-zero scalar $c$ such that $m=m_1+m_2-1$, $(A_1,{\frac{1}{c}}B_1)$ is  $m_1$-symmetric and $(A_2, c B_2)$ is  $m_2$-symmetric \cite[Theorem 5.2]{PG}.

\

In this paper, we start by (equivalently) defining $m$-isometric and $m$-symmetric pairs $(\A,\e)$ of commuting $d$-tuples of Banach space operators in terms of the elementary operators of left and right multiplication (see \cite{DK1}, where this is done for single linear operators). Alongwith introducing other relevant notations and terminology, this is done in Section 2. Section 3 considers the relationship between the $m$-isometric properties of $(\A_i,\e_i)$ (similarly, $m$-symmetric properties of $(\A_i,\e_i)$), $1\leq i\leq 2$, and their product $(\A_1\A_2,\e_1\e_2)$.   A necessary and sufficient condition for product pairs $(\A_1\A_2,\e_1\e_2)$  to be strict $m$-isometric (similarly, strict $m$-symmetric) is proved and its relationship with the strictness of   $m$-isometric pairs $(\A_i,\e_i)$ is explained. Section 4, the penultimate section, proves that the results of Paul and Gu \cite{PG} extend to tensor products of commuting $d$-tuples.  The advantage of our defining $m$-isometric (and, similarly, $m$-symmetric) pairs $(\A,\e)$ using the left/right multiplication operators over definition $(1)$ (resp., $(2)$)  lies in the fact that it provides us with a means to exploit familiar arguments used to prove $1$-tuple (i.e., single linear operator) version of these results.  Here it is seen that the invertibility of $S_i$, $1\leq i\leq 2$, is a sufficient condition, a condition guaranteed by the left invertibility of $S_i$ and the right invertibility of $T_i$ ($1\leq i\leq 2$), in \cite[Theorem 5.2]{PG}.

\section {\sfstp Definitions and introductory properties} For $A,B\in\b$, let $L_A$ and $R_B \in B(\b)$ denote respectively the operators
$$
L_A(X)=AX \ {\rm and}\ R_B(X)=XB
$$
of left multiplication by $A$ and right multiplication by $B$. Given commuting $d$-tuples  $\A=(A_1, \cdots, A_d)$ and $\e=(B_1, \cdots, B_d) \in \b^d$,
let  $\L^{\alpha}_{\A}$ and $\R^{\alpha}_{\e}$,
$$
 \alpha=({\alpha_1, \cdots, \alpha_d}), \ |\alpha|=\sum_{i=1}^d{\alpha_i}, \  \alpha_i\geq 0 \ {\rm for \ all} \  1\leq i\leq d,
$$
 be defined by
$$
\L_{\A}^{\alpha}=\Pi_{i=1}^d{L^{\alpha_i}_{A_i}},  \  \R_{\e}^{\alpha}=\Pi_{i=1}^d{R_{B_i}^{\alpha_i}}.
$$
For an operator $X\in\b$, let convolution  ``*"  and  multiplication  ``$\times$"  denote, respectively, the  operations
\begin{eqnarray*}
& & (\L_{\A} * \R_{\e})^j(X)=\left(\sum_{|\alpha|=j}{{\frac{j!}{{\alpha}!}}\L_{\A}^{\alpha}\R_{\e}^{\alpha}}\right)(X)=\left(\sum_{i=1}^d{L_{A_i}R_{B_i}}\right)^j(X)\\
& & (\rm{all \ integers} \ j\geq 0, \ {\alpha}!={\alpha_1}! \cdots {\alpha_d}!)  \ {\rm and} \\
& & \left(\L_{\A}\times\R_{\e}\right)(X)=\left(\sum_{i=1}^d L_{A_i}\right)\left(\sum_{i=1}^d R_{B_i}\right)(X). 
\end{eqnarray*}
Define the operator $\lfloor \d{A_i X B_i}\rfloor^n$ by 
$$
\lfloor \d{A_i X B_i}\rfloor^n=\d{A_i\lfloor\d {A_i X B_i}\rfloor ^{n-1}B_i} \ {\rm{for \ all \ positive \ integers}} \ n. 
$$
(Thus, $\lfloor A\rfloor^n=A\lfloor A\rfloor^{n-1}I=I\lfloor A\rfloor^{n-1}A= \cdots =A^n$,  $\lfloor  AB\rfloor^n=A\lfloor AB\rfloor^{n-1}B= \cdots  =A^nB^n$ and $\lfloor\d A_iB_i\rfloor^n= \d A_i\lfloor \d A_iB_i\rfloor^{n-1}B_i$.)

\

We say that the $d$-tuples $\A$ and $\e$ commute, $[\A,\e]=0$, if 
$$
[A_i,B_j]=A_iB_j - B_jA_i=0 \ {\rm  for  \  all} \ 1\leq i, j\leq d.
$$  
Evidently,
$$ 
[\L_{\A},\R_{\e}]=0
$$
and if $[\A,\e]=0$, then 
$$
[\L_{\A},\L_{\e}]=[\R_{\A},\R_{\e}]=0.
$$
 A pair $(\A,\e)$ of commuting $d$-tuples $\A$ and $\e$ is said to be  $m$-isometric, $(\A,\e)\in  m$-isometric, for some positive integer $m$, if
\begin{eqnarray*} \triangle^m_{\A,\e}(I) &=& (I-\L_{\A}  *  \R_{\e})^m(I)\\
&=&   \sum_{j=0}^m(-1)^j\left(\begin{array}{clcr}m\\j\end{array}\right)\left(\L_{\A}  *  \R_{\e}\right)^j(I)\\
&=&   \sum_{j=0}^m(-1)^j\left(\begin{array}{clcr}m\\j\end{array}\right)\left(\sum_{i=1}^d L_{A_i}R_{B_i}\right)^j(I)\\
&=&  \sum_{j=0}^m(-1)^j\left(\begin{array}{clcr}m\\j\end{array}\right)\lfloor \d A_iB_i\rfloor^j\\
&=&   \sum_{j=0}^m(-1)^j\left(\begin{array}{clcr}m\\j\end{array}\right)\left(\sum_{|\alpha|=j}{\frac{j!}{{\alpha}!}}\A^{\alpha} \e^{\alpha}\right)\\
&=& 0;
\end{eqnarray*}
$(\A,\e)$ is $n$-symmetric, for some positive integer $n$, if
\begin{eqnarray*} 
\delta_{\A,\e}^n(X) &=& (\L_{\A}-\R_{\e})^n(I)\\
&=& \left(\sum_{j=0}^n(-1)^{j}\left(\begin{array}{clcr}n\\j\end{array}\right) \L_{\A}^{n-j}\times\R_{\e}^j\right)(I)\\
&=& \left(\sum_{j=0}^n(-1)^{j}\left(\begin{array}{clcr}n\\j\end{array}\right)\left(\sum_{i=1}^d L_{A_i}\right)^{n-j}\left(\sum_{i=1}^d R_{B_i}\right)^j\right)(I)\\
&=& \sum_{j=0}^n(-1)^{n-j}\left(\begin{array}{clcr}n\\j\end{array}\right)\left(\sum_{i=1}^d A_i\right)^{n-j}\left(\sum_{i=1}^d B_i\right)^j\\
&=& 0.
\end{eqnarray*}
Commuting tuples of $m$-isometric, similarly $n$-symmetric operators, share a large number of properties with their single operator counterparts: for example, $\triangle^m_{\A,\e}(I)=0$ implies $\triangle^t_{\A,\e}(I)=0$, similarly $\delta^m_{\A,\e}(I)=0$ implies $\delta^t_{\A,\e}(I)=0$, for integers $t\geq m$ 
 However, there are instances where a property holds for the single operator version but fails for the $d$-tuple version: for example, whereas 
$$
 \triangle^m_{A,B}(I)=0\Longrightarrow \triangle^m_{A^t,B^t}(I) \ {\rm for \ all \ integers} \ t\geq 1
$$
 and 
$$
\triangle^m_{A,B}(I)=0\Longleftrightarrow \triangle^m_{A^{-1},B^{-1}}(I)=0 \ {\rm for \ all \ invertible} \ A \ {\rm and} \ B
$$
(similarly, for $m$-symmetric $(A,B)$), these properties fail for $d$-tuples (see \cite{DK1}).

\

If $(\A,\e)\in (X,m)$-isometric, then 
\begin{eqnarray*}
& & \triangle^m_{\A,\e}(X)=0\Longleftrightarrow (I-\L_{\A}  *  R_{\e})\left(\triangle^{m-1}_{\A,\e}(X)\right)=0\\
&\Longleftrightarrow& (\L_{\A}  *  \R_{\e})\triangle^{m-1}_{\A,\e}(X)=\triangle^{m-1}_{\A,\e}(X)\\
&\Longrightarrow& \cdots \Longrightarrow (\L_{\A} * \R_{\e})^t \triangle^{m-1}_{\A,\e}(X)=\triangle^{m-1}_{\A,\e}(X)
\end{eqnarray*}  
and  if $(\A,\e)\in (X,n)$-symmetric, then
\begin{eqnarray*} & &  \delta^n_{\A,\e}(X)=0 \Longleftrightarrow (\L_{\A}-\R_{\e})(\delta^{n-1}_{\A,\e}(X)=0\\
&\Longleftrightarrow& \L_{\A}\delta^{n-1}_{\A,\e}(X)=\R_{\e}(\delta^{n-1}_{\A,\e}(X))\\
&\Longrightarrow& \cdots \Longrightarrow \L^t_{\A}\delta^{n-1}_{\A,\e}(X)=\R^t_{\e}\delta^{n-1}_{\A,\e}(X).
\end{eqnarray*}
for all integers $t\geq 0$. Here 
\begin{eqnarray*}
\L_{\A}(\delta^{n-1}_{\A,\e}(X))&=&\L_{\A}\left(\sum_{j=0}^{n-1}(-1)^{j}\left(\begin{array}{clcr}n-1\\j\end{array}\right) \L_{\A}^{n-1-j}\times\R_{\e}^j\right)(X)\\
&=& \left(\sum_{j=0}^{n-1}(-1)^{j}\left(\begin{array}{clcr}n-1\\j\end{array}\right) \L_{\A}^{n-j}\times\R_{\e}^j\right)(X)
\end{eqnarray*}
and 
\begin{eqnarray*}
\R_{\e}(\delta^{n-1}_{\A,\e}(X)) =\left(\sum_{j=0}^{n-1}(-1)^{j}\left(\begin{array}{clcr}n-1\\j\end{array}\right) \L_{\A}^{n-1-j}\times\R_{\e}^{j+1}\right)(X).
\end{eqnarray*}

\section{\sfstp Results: strictness of products} Let $\A. \e\in \b^d$ be commuting $d$-tuples, and let ${\E}_{\A,\e}$ denote the operator
$$
{\E}_{\A,\e}(X)= (\L_{\A} * \R_{\e})(X), \ X\in\b.
$$
By definition, $(\A,\e)$ is strict $m$-isometric if $\triangle^m_{\A,\e}(I)=0$ and $\triangle^{m-1}_{\A,\e}(I)\neq 0$; similarly, $(\A,\e)$ is strict $m$-symmetric if $\delta^m_{\A,\e}(I)=0$ and $\delta^{m-1}_{\A,\e}(I)\neq 0$. In the following, we give a necessary and sufficient condition for the products pair $(\A_1\A_2, \e_1\e_2)$, $\A_i$ and $\e_i$ commuting $d$-tuples, to be strict $m$-isometric (resp., strict $m$-symmetric), and explore its relationship to the strict $m$-isometric (resp., $m_i$-symmetric) property of $(\A_i,\e_i)$; \  $i=1,2$. We start with a technical lemma.

 \begin{lem}\label{lem00} (i). If $(\A,\e)$ is strict $m$-isometric, then the sequence 
$\{{\E}^{t\pm r}_{\A,\e}\triangle^r_{\A,\e}(I)\}_{r=0}^{m-1}$ is linearly independent for all $t\geq m-1$.

\

\noindent (ii). If $(\A,\e)$ is strict $m$-symmetric and $\L_{\A}^{m-1}\delta^{m-1}_{\A,\e}(I)\neq 0$, then the sequences $\{\L^r_{\A}\delta^r_{\A,\e}(I)\}_{r=0}^{m-1}$ and  $\{\R^r_{\e}\delta^r_{\A,\e}(I)\}_{r=0}^{m-1}$ are linearly independent.
\end{lem}
\begin{demo} The proof in both the cases is by contradiction.

(i). Assume that there exist scalars $a_i$, $0\leq i\leq m-1$, not all zero such that $\sum_{r=0}^{m-1}{a_r{\E}^{t\pm r}_{\A,\e}\triangle^r_{\A,\e}(I)=0}$. Then, since $\triangle^m_{\A,\e}(I)=0$ and ${\E}_{\A,\e}$ commutes with $\triangle_{\A,\e}$,
\begin{eqnarray*}
& & \triangle^{m-1}_{\A,\e}\left(\sum_{r=0}^{m-1}{a_r{\E}^{t\pm r}_{\A,\e}\triangle^r_{\A,\e}(I)}\right)=0\\
&\Longrightarrow& a_0 {\E}^t_{\A,\e}\triangle^{m-1}_{\A,\e}(I)=0\\
&\Longrightarrow& a_0=0,
\end{eqnarray*}
since 
$$
\triangle^m_{\A,\e}(I)=0\Longleftrightarrow \triangle^{m-1}_{\A,\e}(I)={\E}_{\A,\e}\triangle^{m-1}_{\A,\e}(I)
$$ implies 
$$
{\E}_{\A,\e}^t\triangle^{m-1}_{\A\e}(I)=\triangle^{m-1}_{\A,\e}(I)\neq 0.$$
Again,
\begin{eqnarray*}
& & \triangle^{m-2}_{\A,\e}\left(\sum_{r=1}^{m-1}{a_r{\E}^{t\pm r}_{\A,\e}\triangle^r_{\A,\e}(I)}\right)=0\\
&\Longrightarrow& a_1 {\E}^{t\pm 1}_{\A,\e}\triangle^{m-1}_{\A,\e}(I)=0\\
&\Longrightarrow& a_1=0,
\end{eqnarray*}
 and hence repeating the argument 
\begin{eqnarray*}
& & \triangle_{\A,\e}\left(\sum_{r=m-2}^{m-1}{A_r{\E}^{t\pm r}_{\A,\e}\triangle^r_{\A,\e}(I)}\right)=0\\
&\Longrightarrow& a_{m-2} {\E}^{t\pm(m-2)}_{\A,\e}\triangle^{m-1}_{\A,\e}(I)=0\\
&\Longrightarrow& a_{m-2}=0 \Longrightarrow a_{m-1}\triangle^{m-1}_{\A,\e}(I)=0\Longleftrightarrow a_{m-1}=0.
\end{eqnarray*}
This is a contradiction.

\

(ii). We prove the linear independence of the sequence $\{\L^r_{\A}\delta^r_{\A,\e}(\I)\}_{r=0}^{m-1}$; since 
$$
\delta^m_{\A,\e}(I)=0\Longleftrightarrow \L_{\A}\delta^{m-1}_{\A,\e}(I) = \R_{\e}\delta^{m-1}_{\A,\e}(I),
$$
the proof for the linear independence of the second sequence follows from that of the first. Suppose there exist scalars $a_i$, not all zero, such that 
$\sum_{r=0}^{m-1}{a_r\L^r_{\A}\delta^r_{\A,\e}(I)}=0$. Then, since $\delta^m_{\A,\e}(I)=0$, $\L_{\A}$ commutes with $\delta_{\A,\e}$ and $\L_{\A}^{m-1}\delta^{m-1}_{\A,\e}(I)\neq 0$,
\begin{eqnarray*}
& & \delta^{m-1}_{\A,\e}\left(\sum_{r=0}^{m-1}{a_r\L^r_{\A}\delta^r_{\A,\e}(I)}\right)=0\\
&\Longrightarrow& a_0\delta^{m-1}_{\A,\e}(I)=0\Longleftrightarrow a_0=0,\\ & & \delta^{m-2}_{\A,\e}\left(\sum_{r=1}^{m-1}{a_r\L^r_{\A}\delta^r_{\A,\e}(I)}\right)=0\\
&\Longrightarrow& a_1\L_{\A}\delta^{m-1}_{\A,\e}(I)=0\Longleftrightarrow a_1=0,
\end{eqnarray*}
and hence repeating the argument
\begin{eqnarray*}
& & \delta_{\A,\e}\left(\sum_{r=m-2}^{m-1}{a_r\L^r_{\A}\delta^r_{\A,\e}(I)}\right)=0\\
&\Longrightarrow& a_{m-2}\L^{m-2}_{\A}\delta^{m-1}_{\A,\e}(I)=0\Longleftrightarrow a_{m-2}=0\\
&\Longrightarrow& a_{m-1}\L^{m-1}_{\A}\delta^{m-1}_{\A,\e}(I)=0\\
&\Longleftrightarrow& a_{m-1}=0.
\end{eqnarray*}
This is a contradiction.
\end{demo}
\begin{rema}\label{rema00} {\em Since $\delta^{m-1}_{\A,\e}(I)\neq 0$ for a strict $m$-symmetric commuting $d$-tuple $(\A,\e)$, the left invertibility of $\A$ (and hence $\L_{A}$) is a sufficient condition for $\L_{\A}^{m-1}\delta^{m-1}_{\A,\e}(I)\neq 0$.}
\end{rema}
Let $\X\overline{\otimes}\X$ denote the completion, endowed with a reasonable cross norm, of the algebraic tensor product of $\X$ with itself. Let $S\otimes T$ denote the tensor product of $S\in\b$ with $T\in\b$. The tensor product of the $d$-tuples $\A=(A_1, \cdots, A_d)$  and $\e=(B_1, \cdots, B_d)$ is the $d^2$-tuple
$$
\A\otimes\e=(A_1\otimes B_1, \cdots, A_1\otimes B_d,A_2\otimes B_1, \cdots, A_2\otimes B_d, \cdots, A_d\otimes B_1, \cdots, A_d\otimes B_d).
$$
Let $\I=I\otimes I$. (Recall that the operator ${\E}_{\A,\e}$  is defined by ${\E}_{\A,\e}(X)=(\L_{\A} * \R_{\e})(X)$.)

\begin{thm}\label{thm00} Given commuting $d$-tuples ${\A}_i,  {\e}_i\in\b^d$, any two of the following conditions implies the third.

\

(i) $({\A}_1\otimes{\A}_2,{\e}_1\otimes{\e}_2)$ is $m$-isometric; $m=m_1+m_2-1$.

(ii) $({\A}_1,{\e}_1)$ is $m_1$-isometric.

(iii) $({\A}_2,{\e}_2)$ is $m_2$-isometric.
\end{thm}
\begin{demo} It is well known, see \cite{DK1}, that $(ii)$ and $(iii)$ imply $(i)$. We prove $(i)$ and $(ii)$ imply $(iii)$; the proof of $(i)$ and $(iii)$ imply $(ii)$ is similar.

Start by observing that if we let 
$$
{\S}_1={\A}_1\otimes {I}, {\S_2}={I}\otimes{\A}_2, {\T}_1={\e}_1\otimes{I} \ {\rm{and}} \ {\T}_2={I}\otimes{\e}_2,
$$
then 
\begin{eqnarray*}
& & [\S_1,\S_2]=[\T_1,\T_2]=0=[\S_1,\T_2]=[\S_2,\T_1]\\
& & ({\A}_i,{\e}_i)\in m_i-{\rm isometric}\Longleftrightarrow ({\S}_i,{\T}_i)\in m_i-{\rm isometric}, i=1,2,\\
& & ({\A}_1\otimes{\A}_2,{\e}_1\otimes{\e}_2)\in m-{\rm isometric}\Longleftrightarrow ({\S}_1{\S}_2,{\T}_1{\T}_2)\in m-{\rm isometric}
\end{eqnarray*}
and
\begin{eqnarray*}
& &(i)\wedge (ii)\Longrightarrow (iii) \ {\rm {if \ and \ only \ if}} \ ({\S}_1{\S}_2,{\T}_1{\T}_2)\in m-{\rm isometric}\\
& & {\rm and}  \ ({\S}_1,{\T}_1)\in m_1-{\rm {isometric \ imply}} \  ({\S}_2,{\T}_2)\in m_2-{\rm isometric}.
\end{eqnarray*}
Let $t\leq m_1$ be the least positive integer such that $({\S}_1,{\T}_1)\in t$-isometric. (Thus,  $({\S}_1,{\T}_1)$ is strict $t$-isometric.) Then
\begin{eqnarray*}
& &\triangle^m_{{\S}_1{\S}_2,{\T}_1{\T}_2}({{\I}}) = ({I}-{\E}_{{\S}_1{\S}_2,{\T}_1{\T}_2}})^m({{\I})\\
&=& \left({\E}_{{\S}_1,{\T}_1}\triangle_{{\S}_2,{\T}_2}+\triangle_{{\S}_1,{\T}_1})^m\right)({{\I}})\\
&=& \j \left(\begin{array}{clcr}m\\j\end{array}\right){\E}^{m-j}_{{\S}_1,{\T}_1}\triangle^{m-j}_{{\S}_2,{\T}_2}\triangle^j_{{\S}_1,{\T}_1}({{\I}})= 0.
\end{eqnarray*}
Since the operators ${\E}_{{\S}_1,{\T}_1}$ and $\triangle_{{\S}_2,{\T}_2}$ commute, as also do the operators $\triangle_{{\S}_1,{\T}_1}$ and $\triangle_{{\S}_2,{\T}_2}$ , 
\begin{eqnarray*}
& & {\E}^{m-j}_{{\S}_1,{\T}_1}(\triangle^j_{{\S}_1,{\T}_1})(\I)=(-1)^j {\E}^{m-j}_{\S_1,\T_1}({\E}_{\S_1,\T_1}-I)^j(\I)\\
&=& \left(\sum_{k=0}^j (-1)^{j+k}\left(\begin{array}{clcr}j\\k\end{array}\right){\E}^{m-k}_{{\S}_1,{\T}_1}\right)({\I})\\
&=& \sum_{k=0}^j (-1)^{j+k}\left(\begin{array}{clcr}j\\k\end{array}\right)\lfloor \d{L_{A_{1i}\otimes I}R_{{B_{1i}\otimes I}}}\rfloor^{m-k}(\I)\\
&=& \sum_{k=0}^j (-1)^{j+k}\left(\begin{array}{clcr}j\\k\end{array}\right)\lfloor\d A_{1i}B_{1i}\otimes I\rfloor^{m-k}\\
& & (\lfloor\d A_{1i}B_{1i}\otimes I \rfloor^n=\d (A_{1i}\otimes I) \lfloor\d A_{1i}B_{1i}\otimes I\rfloor^{n-1}(B_{1i}\otimes I),\\
& & {\rm all \ positive \ integers} \ n)\\
&=& \sum_{k=0}^j (-1)^{j+k}\left(\begin{array}{clcr}j\\k\end{array}\right)\lfloor\d A_{1i}B_{1i}\rfloor^{m-k}\otimes I\\
&=& X_j\otimes I \ ({\rm say}).
\end{eqnarray*}
Hence
\begin{eqnarray*} 
& & \triangle^{m-j}_{{\S}_2,{\T}_2}\left( {\E}^{m-j}_{{\S}_1,{\T}_1}(\triangle^j_{{\S}_1,{\T}_1}(\I)\right)=(-1)^{m-j}({\E}_{\S_2,\T_2}-I)^{m-j}\left( {\E}^{m-j}_{{\S}_1,{\T}_1}\triangle^j_{{\S}_1,{\T}_1}(\I)\right)\\
&=& \sum_{p=0}^{m-j} (-1)^{m-j+p}\left(\begin{array}{clcr}m-j\\p\end{array}\right)\left(I\otimes\lfloor\d A_{2i}B_{2i}\rfloor\right)^{m-j-p}\left(X_j\otimes I\right)\\
&=& \sum_{p=0}^{m-j} (-1)^{m-j+p}\left(\begin{array}{clcr}m-j\\p\end{array}\right)\left(X_j\otimes\lfloor\d A_{2i}B_{2i}\rfloor^{m-j-p}\right)\\
&=& X_j\otimes Y_j \ {\rm say},
\end{eqnarray*}
and
$$ \triangle^m_{{\S}_1{\S}_2,{\T}_1{\T}_2}(\I)= \sum_{j=0}^{m} \left(\begin{array}{clcr}m\\j\end{array}\right) (X_j\otimes Y_j).
$$
The sequence $\{X_j\}_{j=0}^{t-1}$ being linearly independent, we must have $Y_j=0$ for all $0\leq j\leq t-1$. Since $j\leq t-1$ implies $m-j=m_1+m_2-1-t+1\geq m_2$, we have 
\begin{eqnarray*}
& & I\otimes \sum^{m_2}_{p=0} (-1)^{m_2+p}\left(\begin{array}{clcr}m_2\\p\end{array}\right) \lfloor\d{A_{2i}B_{2i}}\rfloor^{m_2-p}=0\\
&\Longleftrightarrow& \sum^{m_2}_{p=0} (-1)^{m_2+p}\left(\begin{array}{clcr}m_2\\p\end{array}\right) \lfloor\d{A_{2i}B_{2i}}\rfloor^{m_2-p}=\triangle^{m_2}_{{\A}_2,{\e}_2}(I)=0.
\end{eqnarray*}
This completes the proof.
\end{demo}

An anlogue of Theorem \ref{thm00} holds for products of $m$-symmetric operators.
\begin{thm}\label{thm10} Given commuting $d$-tuples $\A_i, \e_i\in\b^d$ such that $\A_i$ is left invertible for $1\leq i\leq 2$, any two of the following conditions implies the third.

\

(i). $(\A_1\otimes\A_2,\e_1\otimes\e_2)$ is $m$-symmetric; $m=m_1+m_2-1$.

(ii). $(\A_1,\e_1)$ is $m_1$-symmetric.

(iii). $(\A_2,\e_2)$ is $m_2$-symmetric.
\end{thm}
\begin{demo} That $(ii)$ and $(iii)$ imply $(i)$, without any hypothesis on the left invertibility of $\A_1$ and $\A_2$, is well known \cite{DK1}. We prove $(i)$ and $(ii)$ imply $(iii)$; the proof of $(i)$ and $(iii)$ imply $(ii)$ is similar and left to the reader.

\

Assume $t\leq m_1$ is the least positive integer such that $\delta^{t-1}_{\A_1,\e_1}({I})\neq 0$. (Thus, ${\A_1,\e_1}$ is strict $t$-symmetric.)  Then
\begin{eqnarray*}
\delta^m_{\A_1\otimes\A_2,\e_1\otimes\e_2}(\I) &=& (\L_{\A_1\otimes\A_2} -\R_{\e_1\otimes\e_2})^m(\I)\\
&=& \left(\L_{\A_1}\otimes\L_{\A_2} - \R_{\e_1}\otimes\R_{\e_2}\right)^m(\I)\\
&=& \left( (\L_{\A_1}\otimes\L_{\A_2} - \R_{\e_1}\otimes\L_{\A_2}) + (\R_{\e_1}\otimes\L_{\A_2} - \R_{\e_1}\otimes\R_{\e_2})\right)^m(\I)\\
&=& \left( \sum_{j=0}^m \left(\begin{array}{clcr}m\\j\end{array}\right) \left(\delta^{m-j}_{\A_1,\e_1}\otimes \L^{m-j}_{\A_2}\right) \times \left(\R^j_{\e_1}\otimes \delta^j_{\A_2,\e_2}\right)\right)(\I)\\
&=& \j \left(\begin{array}{clcr}m\\j\end{array}\right) \left(\R^j_{\e_1} \delta^{m-j}_{\A_1,\e_1}(I)\right) \otimes \left(\L^{m-j}_{\A_2} \delta^j_{\A_2,\e_2}(I)\right).
\end{eqnarray*}
Since $(\A_1,\e_1)$ is strict $t$-symmetric and $\L_{\A_1}\delta^{t-1}_{\A_1,\e_1}(I)=\R_{\e_1}\delta^{t-1}_{\A_1,\e_1}(I)\neq 0$, the argument of the proof of Lemma \ref{lem00} implies the linear independence of the sequence $\{\R^j_{\e_1}\delta^{m-j}_{\A_1,\e_1}(I)\}_{m-j=0}^{t-1}$. Hence 
$\L^{m-j}_{\A_2}\delta^j_{\A_2,\e_2}(I)=0$ for all $m-j\leq t-1$, equivalently, $j\geq m-t+1\geq m_1+m_2-1-m_1+1=m_2$. But then, since $\L_{\A_2}$ is left invertible, $\delta^j_{\A_2,\e_2}(I)=0$ for all $j\geq m_2$.
\end{demo}

Strictness in conditions $(i)  -  (iii)$ of Theorem \ref{thm00} requires more: thus whereas $(\A_1\otimes\A_2,\e_1\otimes\e_2)$ is strictly $m$ isometric implies $(\A_i,\e_i)$ is strictly $m_i$-isometric for both $i=1$ and $i=2$, $(\A_i,\e_i)$ is strictly $m_i$-isometric for both $i=1$ and $i=2$ does not in general imply
$(\A_1\otimes\A_2,\e_1\otimes\e_2)$ is strictly $m$ isometric. 
\begin{thm}\label{thm01} Given commuting $d$-tuples $\A_i, \e_i\in\b^d$, $1\leq i\leq 2$, such that $[\A_1,\A_2]=[\e_1,\e_2]=0$, if $(\A_i,\e_i)$ is  $m_i$-isometric and $(\A_1\A_2,\e_1\e_2)$ is  $m$ isometric, $m=m_1+m_2-1$, then:

\

\noindent (i) $(\A_1\A_2,\e_1\e_2)$ is strictly $m$ isometric if and only if
\begin{eqnarray}
\triangle^{m_1-1}_{\A_1,\e_1}\left(\triangle^{m_2-1}_{\A_2,\e_2}(I)\right)= \triangle^{m_2-1}_{\A_2,\e_2}\left(\triangle^{m_1-1}_{\A_1,\e_1}(I)\right)\neq 0;
\end{eqnarray}
(ii) $(\A_1\A_2,\e_1\e_2)$ is strictly $m$ isometric implies $(\A_i,\e_i)$ is strictly isometric for $1\leq i\leq 2$;

\noindent (iii) $(\A_i,\e_i)$ is strictly isometric for $1\leq i\leq 2$ does not imply $(\A_1\A_2,\e_1\e_2)$ is strictly $m$ isometric
\end{thm}
\begin{demo} The hypothesis $[\A_1,\A_2]=[\e_1,\e_2]=0$ implies
$$
[\triangle_{\A_1,\e_1},\triangle_{\A_2,\e_2}]=0=[{\E}^{m_2-1}_{\A_1,\e_1},\triangle_{\A_2,\e_2}].
$$
Since
\begin{eqnarray*} & & {{\E}^{m_2-1}_{\A_1,\e_1}}\triangle^{m_1-1}_{\A_1,\e_1}\left(\triangle^{m_2-1}_{\A_2,\e_2}(I)\right)\\
&=&  \triangle^{m_2-1}_{\A_2,\e_2}\left({\E}^{m_2-1}_{\A_1,\e_1}\triangle^{m_1-1}_{\A_1,\e_1}(I)\right)\\
&=& \triangle^{m_2-1}_{\A_2,\e_2}\left(\triangle^{m_1-1}_{\A_1,\e_1}(I)\right),
\end{eqnarray*}
whenever $\triangle^{m_1}_{\A_1,\e_1}(I)=0$ (equivalently, ${\E}_{\A_1,\e_1} \triangle^{m_1-1}_{\A_1,\e_1}(I)= \triangle^{m_1-1}_{\A_1,\e_1}(I)$), if $(\A_1\A_2,\e_1\e_2)$ is $m$-isometric and either of $(\A_i,\e_i)$, $1\leq i\leq 2$, is not strict $m_i$-isometric then $(3)$ is contradicted. Hence $(i)$ implies $(ii)$.

\

To prove $(i)$, assume $(\A_1\A_2,\e_1\e_2)$ is strict $m$-isometric . Then $\triangle^m_{\A_1\A_2,\e_1\e_2}(I)=0$ and $\triangle^{m-1}_{\A_1\A_2,\e_1\e_2}(I)\neq 0$. We have
\begin{eqnarray*}
0 &\neq & \triangle^{m-1}_{\A_1\A_2,\e_1\e_2}(I)\\
&=& \sum_{j=0}^{m-1}  \left(\begin{array}{clcr}m-1\\j\end{array}\right) \triangle^{m-1-j}_{\A_2,\e_2}\left( {\E}^{m-1-j}_{\A_1,\e_1}\triangle^j_{\A_1,\e_1}(I)\right)\\
& & ({\rm see \ the \ proof \ of Theorem \ \ref{thm00}})\\
&=& \sum_{j=0}^{m-1} \left(\begin{array}{clcr}m-1\\j\end{array}\right)      \triangle^{m-1-j}_{\A_2,\e_2}\left({\E}^{m-1-j}_{\A_1,\e_1}\triangle^j_{\A_1,\e_1}(I)\right)\\
& & ({\rm since} \ \triangle^{m_1}_{\A_1,\e_1}(I)=0 \ {\rm for} \ j\geq m_1)\\
&=& \left(\begin{array}{clcr}m-1\\m_1-1\end{array}\right) \triangle^{m_2-1}_{\A_2,\e_2}\left(\triangle^{m_1-1}_{\A_1,\e_1}(I)\right),
\end{eqnarray*}
since $\triangle^{m-1-j}_{\A_2,\e_2}(I)=0$ for all $m-1-j\geq m_2$, equivalently, $m_1-2\geq j$, and ${\E}^{m-1-(m_1-1)}_{\A_1,\e_1}\left(\triangle^{m_1-1}_{\A_1,\e_1}(I)\right)=\triangle^{m_1-1}_{\A_1,\e_1}(I)$.

\

To complete the proof, we give an example proving $(iii)$. Let $I_2=I\oplus I$, and let $A, B\in\b$ be such that $(A,B)$ is strict $m$-isometric. Define operators $\A_i, \e_i\in B(\X\oplus\X)^d$, $1\leq i\leq 2$, by
\begin{eqnarray*}
& & \A_1 = (A_{11}, \cdots, A_{1d}) = {\frac{1}{\sqrt{d}}}(A\oplus I, \cdots, A\oplus I),\\
& & \A_2 = (A_{21}, \cdots, A_{2d}) = {\frac{1}{\sqrt{d}}}(I\oplus A, \cdots, I\oplus A),\\
& & \e_1 = (B_{11}, \cdots, B_{1d}) = {\frac{1}{\sqrt{d}}}(B\oplus I,  \cdots, B\oplus I),\\
& & \e_2 = (B_{21}, \cdots, B_{2d}) = {\frac{1}{\sqrt{d}}}(I\oplus B, \cdots, I\oplus B).
\end{eqnarray*}
Then 
\begin{eqnarray*}
\triangle^m_{\A_1,\e_1}(I_2) &=&  \sum_{j=0}^m (-1)^j\left(\begin{array}{clcr}m\\j\end{array}\right)\left(\L_{\A_1} * \R_{\e_1}\right)^{m-j}(I_2)\\
&=& \sum_{j=0}^m (-1)^j\left(\begin{array}{clcr}m\\j\end{array}\right)\lfloor \d {\frac{1}{d}}(AB\oplus I)\rfloor^{m-j}\\
& & \left(\lfloor\d (AB\oplus I)\rfloor^t=\d (A\oplus I)\lfloor\d (AB\oplus I)\rfloor^{t-1}(B\oplus I)  \ {\rm {for \ integers}} \ t\geq 1\right)\\
&=& \sum_{j=0}^m (-1)^j\left(\begin{array}{clcr}m\\j\end{array}\right)\lfloor (AB\oplus I)\rfloor^{m-j}\\
&=& \sum_{j=0}^m (-1)^j\left(\begin{array}{clcr}m\\j\end{array}\right) (A^{m-j}B^{m-j}\oplus I)\\
&=& 0,
\end{eqnarray*}
and 
$$ 
\triangle^{m-1}_{\A_1,\e_1}(I_2)\neq 0.
$$
Similarly, $\triangle^{m}_{\A_2,\e_2}(I_2)=0$  and $\triangle^{m-1}_{\A_2,\e_2}(I_2)\neq 0$.

\

Consider now $\triangle^m_{\A_1\A_2,\e_1\e_2}(I_2)$. Since $\A_i$ and $\e_i$ are commuting $d$-tuples such that $[\A_1,\A_2]=[\e_1,\e_2]=0$, $(\A_1\A_2,\e_1\e_2)$ is  $(2m-1)$-isometric. Again, since 
$$
\A_1\A_2={\frac{1}{d}} (A\oplus A, \cdots, A\oplus A)\in B(\X\oplus\X)^{d^2}
$$
and
$$
\e_1\e_2={\frac{1}{d}} (B\oplus B, \cdots, B\oplus B)\in B(\X\oplus\X)^{d^2},
$$
 \begin{eqnarray*}
& &\triangle^m_{\A_1\A_2,\e_1\e_2}(I_2)\\
 &=&  \sum_{j=0}^m (-1)^j\left(\begin{array}{clcr}m\\j\end{array}\right)\left(\L_{\A_1\A_2} * \R_{\e_1\e_2}\right)^{m-j}(I_2)\\
&=& \sum_{j=0}^m (-1)^j\left(\begin{array}{clcr}m\\j\end{array}\right)  \lfloor\sum_{i=1}^{d^2}{\frac{1}{d^2}} (AB\oplus AB)\rfloor^{m-j}\\
& & (\lfloor\sum_{i=1}^{d^2}(AB\oplus AB)\rfloor^t=\sum_{i=1}^{d^2}(A\oplus A)\lfloor\sum_{i=1}^{d^2}(AB\oplus AB)\rfloor^{t-1}(B\oplus B) \\
& & {\rm for \ all \ integers }\ t\geq 1)\\
&=&   \sum_{j=0}^m (-1)^j\left(\begin{array}{clcr}m\\j\end{array}\right) \lfloor (AB\oplus AB)\rfloor^{m-j}\\
&=&  \sum_{j=0}^m (-1)^j\left(\begin{array}{clcr}m\\j\end{array}\right) (A^{m-j}B^{m-j}\oplus A^{m-j}B^{m-j})\\
&=& 0,
\end{eqnarray*}
i.e., $(\A_1\A_2,\e_1\e_2)$ is $m$-isometric. Thus $(\A_1\A_2,\e_1\e_2)$ is not strict $(2m-1)$-isometric for all $m\geq 2$.
\end{demo}
\begin{rema}\label{rema01} {\em Choosing $A, B\in\b$, and $\A_i, \e_i\in B(\X\oplus\X)^d$, to be the operators of the example proving part $(iii)$, define operators $\S_i$ and $\T_i \in B((\X\oplus\X)\overline{\otimes}(\X\oplus\X))^d$ by $\S_i=\A_i\otimes I_2$ and $\T_i=\e_i\otimes I_2$; $1\leq i\leq 2$ and $I_2=I\oplus I$.
Then $(\S_i,\T_i)$ and $(\S_1\otimes\S_2,\T_1\otimes\T_2)$, $1\leq i\leq 2$, are all strict $m$-isometric.}
\end{rema} 
Just as for $m$-isometric operators, strictness for $m$-symmetric operators requires more.
\begin{thm}\label{thm11}  Given commuting $d$-tuples $\A_i, \e_i\in \b^d$ such that $[\A_1,\A_2]=[\e_1,\e_2]=0$, and $\A_i$ is left invertible for $1\leq i\leq 2$, if $(\A_i,\e_i)$ is $m_i$-symmetric, $1\leq i\leq 2$, and $(\A_1\A_2,\e_1\e_2)$ is $m$-symmetric, $m=m_1+m_2-1$, then:

\

\noindent (i) $(\A_1\A_2,\e_1\e_2)$ is strict $m$-symmetric if and only if
\begin{eqnarray}
& &  \L_{\A_1} ^{m_2-1}\times \delta^{m_2-1}_{\A_2,\e_2}\left(\R^{m_1-1}_{\A_2} \times \delta^{m_1-1}_{\A_1,\e_1}(I)\right)\\&=& \R^{m_1-1}_{\A_2} \times \delta^{m_1-1}_{\A_1,\e_1}\left(\L_{\A_1} \times \delta^{m_2-1}_{\A_2,\e_2}(I)\right)\neq 0;
\end{eqnarray}
(ii) $(\A_1\A_2,\e_1\e_2)$ is strict $m$-symmetric implies $(\A_i,\e_i)$ strict $m_i$-symmetric for $1\leq i\leq 2$;

\noindent (iii) $(\A_i,\e_i)$ is strict $m_i$-symmetric for $1\leq i\leq 2$ does not imply $(\A_1\A_2,\e_1\e_2)$ is strict $m$-symmetric.
\end{thm}
\begin{demo}  $(i)$ and $(ii)$. Evidently, if either of $(\A_i,\e_i)$ is not strict $m_i$-symmetric, then $(4)$ and $(5)$  equal $0$ and  $(i)$ is violated. Hence $(i)$ implies $(ii)$. We prove $(i)$, and then modify the example in the proof of Theorem \ref{thm01} to prove $(iii)$.

 $(\A_1\A_2,\e_1\e_2)$ is strict $m$-symmetric if and only if $\delta^m_{\A_1\A_2,\e_1\e_2}(I)=0$ and 
\begin{eqnarray*}
0 &\neq& \delta^{m-1}_{\A_1\A_2,\e_1\e_2}(I)=(\L_{\A_1\A_2} - \R_{\e_1\e_2})^{m-1}(I)\\
&=&  \sum_{j=0}^{m-1} \left(\begin{array}{clcr}m-1\\j\end{array}\right) \left(L^{m-1-j}_{\A_1} \times \delta^{m-1-j}_{\A_2,\e_2} \right)  \left(\R^j_{\A_2} \times \delta^j_{\A_1,\e_1}\right)(I)\\
&=&  \sum_{j=0}^{m-1} \left(\begin{array}{clcr}m-1\\j\end{array}\right)  \left(\R^j_{\A_2} \times \delta^j_{\A_1,\e_1}\right)\left(L^{m-1-j}_{\A_1} \times \delta^{m-1-j}_{\A_2,\e_2} \right)(I)
\end{eqnarray*}
by the commutativity hypotheses on $\A_i, \e_i$ and the commutativity of the left and the right multiplication operators. Since
$$
\R^j_{\A_2} \times \delta^j_{\A_1,\e_1}(I)=(\d R_{A_{2i}})^j(\delta^j_{\A_1,\e_1}(I))=0
$$
for all $j\geq m_1$
,
 $$
0\neq  \sum_{j=0}^{m_1-1} \left(\begin{array}{clcr}m-1\\j\end{array}\right)  \left(\R^j_{\A_2} \times \delta^j_{\A_1,\e_1}\right)\left(L^{m-1-j}_{\A_1} \times \delta^{m-1-j}_{\A_2,\e_2} \right)(I).
$$
But then, since $\delta^{m-1-j}_{\A_2,\e_2}(I)=0$ for $m-1-j=m_1+m_2-2-j \geq m_2$,
$$
0\neq \left(\begin{array}{clcr}m-1\\m_1-1\end{array}\right)  \left(\R^{m_1-1}_{\A_2} \times \delta^{m_1-1}_{\A_1,\e_1}\right)\left(L^{m_2-1}_{\A_1} \times \delta^{m_2-1}_{\A_2,\e_2} \right)(I).
$$

\

$(iii)$. Define operators $\A_i$ and $\e_i$, $1\leq i\leq 2$, as in the proof of Theorem \ref{thm01}(iii). Choose the operators $A, B$ this time to be such that $A$ is left invertible and $(A,B)$ is strict $m$-symmetric. Let, as before,  $I_2=I\oplus I$. Then, since
\begin{eqnarray*}
\delta^m_{\A_1,\e_1}(I_2) &=& \left( \sum_{j=0}^{m_1-1} \left(\begin{array}{clcr}m-1\\j\end{array}\right)  \left(\d L_{A_{1i}}\right)^{m-j} \left(\d R_{B_{1i}}\right)^j\right)(I_2)\\
&=& \left( \sum_{j=0}^{m_1-1} \left(\begin{array}{clcr}m-1\\j\end{array}\right) \left(L_{A\oplus I}^{m-j}R^j_{B\oplus I}\right)\right)(I_2)\\
&=& \left( \sum_{j=0}^{m_1-1} \left(\begin{array}{clcr}m-1\\j\end{array}\right) \left(L_A^{m-j}R^j_B \oplus I\right)\right)(I_2)\\
&=& \delta^m_{A,B}(I)\oplus 0=0,
\end{eqnarray*}
i.e., $(\A_1,\e_1)$ is $m$-symmetric. Similarly, $(\A_2,\e_2)$ is $m$-symmetric. This, in view of the fact that $[\A_1,\A_2]=[\e_1,\e_2]=0$, implies $(\A_1\A_2,\e_1\e_2)$ is $(2m-1)$-symmetric. However,
\begin{eqnarray*}
\delta^{m}_{\A_1\A_2,\e_1\e_2}(I_2) &=& \left(\sum_{j=0}^{m-1} \left(\begin{array}{clcr}m-1\\j\end{array}\right) \left({\frac{1}{d^2}}\sum_{i=1}^{d^2}L_{A\oplus A}\right)^{m-j}\left({\frac{1}{d^2}}\sum_{i=1}^{d^2}R_{\e\oplus\e}\right)^j\right)(I_2)\\
&=& \left(\sum_{j=0}^m\left(\begin{array}{clcr}m-1\\j\end{array}\right)\left(L^{m-j}_AR^j_B\oplus L^{m-j}_AR^j_B\right)\right)(I_2)\\
&=& 0.
\end{eqnarray*}
Hence $(\A_1\A_2,\e_1\e_2)$ is not strict $(2m-1)$-symmetric for all $m> 1$.
\end{demo}

\section{\sfstp Results: the inverse problem}By definition, $(S_1\otimes S_2, T_1\otimes T_2)$ is $m$-isometric if and only if
\begin{eqnarray*} \triangle^m_{S_1\otimes S_2, T_1\otimes T_2}(I\otimes I) &= & \left( \sum_{j=0}^{m} (-1)^j\left(\begin{array}{clcr}m\\j\end{array}\right) \left(L_{S_1\otimes S_2} R_{T_1\otimes T_2}\right)^j\right)(I\otimes I)\\
&=&  \sum_{j=0}^{m}(-1)^j \left(\begin{array}{clcr}m\\j\end{array}\right)\lfloor S_1T_1\rfloor^j\otimes \lfloor S_2T_2\rfloor^j\\
&=& 0
\end{eqnarray*}
and it is strict $m$-isometric if and only if it is $m$-isometric and 
$$
 \triangle^{m-1}_{S_1\otimes S_2, T_1\otimes T_2}(I\otimes I)= \sum_{j=0}^{m-1} (-1)^j\left(\begin{array}{clcr}m-1\\j\end{array}\right)\lfloor S_1T_1\rfloor^j\otimes \lfloor S_2T_2\rfloor^j\neq 0.
$$
(Recall that, given operators $S_i, T_i\in\b$, $1\leq i\leq d$,  $\lfloor\d S_iT_i\rfloor^t=S_1\lfloor\d S_iT_i\rfloor^{t-1}T_1+ \cdots +S_d\lfloor\d S_iT_i\rfloor^{t-1}T_d$.) Paul and Gu \cite[Theorem 1.1]{PG} prove that ``if $(S_1\otimes S_2, T_1\otimes T_2)$ is  $m$-isometric, then there exist integers $m_i >0$, and a non-zero scalar $c$, such that $m=m_1+m_2-1$, $(S_1, {\frac{1}{c}}T_1)$ is $m_1$-isometric and $(S_2,cT_2)$ is $m_2$-isometric". Translating this into the terminology 
above, one has the following.
\begin{pro}\label{pro00} Given operators $S_i, T_i\in \b$, $1\leq i\leq 2$, if 
$$
\sum_{j=0}^{m} (-1)^j\left(\begin{array}{clcr}m\\j\end{array}\right)\lfloor S_1T_1\rfloor^j\otimes \lfloor S_2T_2\rfloor^j=0
$$ 
then there exist integers $m_i >0$, and a non-zero scalar $c$, such that $m=m_1+m_2-1$,

$$
 \sum_{j=0}^{m_1} (-1)^j\left(\begin{array}{clcr}m_1\\j\end{array}\right)\lfloor S_1 ({\frac{1}{c}}T_1)\rfloor^j=0=\sum_{j=0}^{m_2} (-1)^j\left(\begin{array}{clcr}m_2\\j\end{array}\right) \lfloor S_2 (cT_2)\rfloor^j.
$$
\end{pro}
The following theorem is an analogue of \cite[Theorem 1.1]{PG} for commuting $d$-tuples of operators. Our proof, which depends on an application of Proposition \ref{pro00}, is achieved by reducing the problem to that for single linear operators.  
\begin{thm}\label{thm12} If $\S_i, \T_i\in\b^d$, $1\leq i\leq 2$, are commuting $d$-tuples such that $(\S_1\otimes\S_2,\T_1\otimes\T_2)$ is strict $m$-isometric, then there exist integers $m_i >0$, and a non-zero scalar $c$, such that $m=m_1+m_2-1$, $(\S_1,{\frac{1}{c}}\T_1)$ is strict $m_1$-isometric and $(\S_2, c\T_2)$ is strict $m_2$-isometric.
\end{thm} 
\begin{demo} If $(\S_1\otimes\S_2,\T_1\otimes\T_2)$ is strict $m$-isometric, then,
since 
$$
\S_1\otimes\S_2=(S_{11}\otimes S_{21}, \cdots, S_{11}\otimes S_{2d},S_{12}\otimes S_{21}, \cdots ,S_{12}\otimes S_{2d},\cdots, S_{1d}\otimes S_{21}, \cdots, S_{1d}\otimes S_{2d})
$$
and
$$
\T_1\otimes\T_2=(T_{11}\otimes T_{21}, \cdots, T_{11}\otimes T_{2d},T_{12}\otimes T_{21}, \cdots, T_{12}\otimes T_{2d}, \cdots, T_{1d}\otimes T_{21}, \cdots, T_{1d}\otimes T_{2d}),
$$ 
\begin{eqnarray*} 0 &=& \triangle^m_{\S_\otimes\S_2,\T_\otimes\T_2}(\I)\\
&=&\left(\sum_{j=0}^{m}(-1)^j \left(\begin{array}{clcr}m\\j\end{array}\right)\left(\L_{\S_1\otimes\S_2} * \R_{\T_1\otimes\T_2}\right)^j\right)(\I\\
&=& \sum_{j=0}^{m}(-1)^j \left(\begin{array}{clcr}m\\j\end{array}\right)\lfloor\sum_{i,k=1}^d(S_{1i}\otimes S_{2k})(T_{1i}\otimes T_{2k})\rfloor^j.
\end{eqnarray*}

This, since
\begin{eqnarray*}
& & \lfloor\sum_{i,k=1}^d(S_{1i}\otimes S_{2k})(T_{1i}\otimes T_{2k})\rfloor=\lfloor\sum_{i,k=1}^d S_{1i}T_{1i}\otimes S_{2i}T_{2k}\rfloor\\
&=& \lfloor\d S_{1i}T_{1i}\rfloor \otimes \lfloor \sum_{k=1}^d S_{2k}T_{2k}\rfloor,
\end{eqnarray*}
implies
\begin{eqnarray*}
0 &=& \triangle^m_{\S_1\otimes\S_,\T_1\otimes\T_2}(\I)\\
&=&  \sum_{j=0}^{m} (-1)^j\left(\begin{array}{clcr}m\\j\end{array}\right)\lfloor\d S_{1i}T_{1i}\rfloor^j \otimes \lfloor \sum_{k=1}^d S_{2k}T_{2k}\rfloor^j\\
\end{eqnarray*}
and
\begin{eqnarray*}
 0\neq & & \triangle^{m-1}_{\S_1\otimes\S_,\T_1\otimes\T_2}(\I)\\
&=&  \sum_{j=0}^{m-1} (-1)^j\left(\begin{array}{clcr}m-1\\j\end{array}\right)\lfloor\d S_{1i}T_{1i}\rfloor^j \otimes \lfloor \sum_{k=1}^d S_{2k}T_{2k}\rfloor^j.
\end{eqnarray*}
Applying Proposition \ref{pro00} we have the existence of a non-zero scalar $c$ and positive integers $m_i$, $m=m_1+m_2-1$, such that 
$$
\sum_{j=0}^{m_1} (-1)^j\left(\begin{array}{clcr}m_1\\j\end{array}\right)\lfloor\d S_{1i}({\frac{1}{c}}T_{1i})\rfloor^j=0=\sum_{j=0}^{m_2}(-1)^j \left(\begin{array}{clcr}m_2\\j\end{array}\right)\lfloor\d S_{2i}(c T_{1i})\rfloor^j.
$$
The strictness of the $m$-isometric property of the tensor products pair $(S_1\otimes S_2, T_1\otimes T_2)$ implies 
$$
 \sum_{j=0}^{m_1-1} (-1)^j\left(\begin{array}{clcr}m_1-1\\j\end{array}\right)\lfloor\d S_{1i}({\frac{1}{c}}T_{1i})\rfloor^j \neq 0\neq \sum_{j=0}^{m_2-1}(-1)^j \left(\begin{array}{clcr}m_2-1\\j\end{array}\right)\lfloor\d S_{2i}(c T_{2i})\rfloor^j,
$$
i.e.,  $(\S_1,{\frac{1}{c}}\T_1)$ is strict $m_1$-isometric and $(\S_2,c\T_2)$ is strict $m_2$-isometric.
\end{demo}

Observe from the proof above that the strictness property of the $m$-isometric operator pair $(\S_1\otimes\S_2,\T_1\otimes\T_2)$ plays no role in the determination of the scalar $c$ or positive integers $m_i$ such $(\S_1,{\frac{1}{c}}\T_1)$ is $m_1$-isometric and $(\S_2,c\T_2)$ is $m_2$-isometric: strictness plays a role only in the determining of the strictness of the $m_1$-isometric property of  $(\S_1,{\frac{1}{c}}\T_1)$  and  the strictness of the  $m_2$-isometric property of $(\S_2,c\T_2)$.

\

If an operator $T\in\B$ is $m$-symmetric, then $\sigma_a(T)$, the approximate point spectrum of $T$, is a subset of $\R$. Hence $T$ is left invertible if and only if it is invertible. Consider operators $S_i, T_i\in\b$, $1\leq i\leq 2$, such that $S_i$ is invertible and $\delta^m_{S_1\otimes S_2,T_1\otimes T_2}(I\otimes I)=0$. Then
\begin{eqnarray*}
& & \delta^m_{S_1\otimes S_2,T_1\otimes T_2}(I\otimes I)= \left(\sum_{j=0}^{m} (-1)^j\left(\begin{array}{clcr}m\\j\end{array}\right)L^{m-j}_{S_1\otimes S_2}R^j_{T_1\otimes T_2}\right)(I\otimes I)=0\\
&\Longleftrightarrow& \left(\sum_{j=0}^{m}(-1)^j \left(\begin{array}{clcr}m\\j\end{array}\right)L^{-j}_{S_1\otimes S_2}R^j_{T_1\otimes T_2}\right)(I\otimes I)=0\\
&\Longleftrightarrow& \triangle^m_{S^{-1}_1\otimes S^{-1}_2,T_1\otimes T_2}(I\otimes I)=0.
\end{eqnarray*}
Assuming, further, $(S_1\otimes S_2,T_1\otimes T_2)$ to be strict $m$-symmetric, it follows that 
$$
(S_1\otimes S_2,T_1\otimes T_2) \  {\rm is \  strict} \  m-{\rm symmetric} \Longleftrightarrow (S^{-1}_1\otimes S^{-1}_2,T_1\otimes T_2) \  {\rm is \  strict} \  m-{\rm isometric}.
$$
Hence there  exists  a non-zero scalar $c$ and positive integers $m_i$, $m=m_1+m_2-1$, such that 
$$
(S^{-1}_1,{\frac{1}{c}}T_1)  \ {\rm is \ strict} \ m_1-{\rm isometric \ and} \ (S_2^{-1},cT_2) \ {\rm is \ strict} \ m_2-{\rm isometric}.
$$
Since
\begin{eqnarray*}
\triangle^m_{S^{-1}_i,\alpha T_i}( I)=0 &\Longleftrightarrow& \left(\sum_{j=0}^{m_i} (-1)^j\left(\begin{array}{clcr}m_i\\j\end{array}\right)L^{-j}_{S_i}(\alpha R_{T_i})^j\right)(I)=0\\
&\Longleftrightarrow&  \left(\sum_{j=0}^{m_i} (-1)^j\left(\begin{array}{clcr}m_i\\j\end{array}\right)L^{m-j}_{S_i}(\alpha R_{T_i})^j\right)(I)=0\\
&\Longleftrightarrow& \delta^{m_i}_{S_i,\alpha T_i}(I)=0
\end{eqnarray*}
and, similarly,
$$
\triangle^{m_i-1}_{S^{-1}_i,\alpha T_i}(I)\neq 0\Longleftrightarrow \delta^{m_i-1}_{S_i,\alpha T_i}(I)\neq 0,
$$
we have the following Banach space analogue of \cite[Theorem 1.2]{PG}.
\begin{pro}\label{pro01} If $S_1, S_2\in\b$ are invertible and $(S_1\otimes S_2,T_1\otimes T_2)$ is  $m$-symmetric, for some operators $T_1, T_2\in\b$, then there  exists  a non-zero scalar $c$ and positive integers $m_i$, $m=m_1+m_2-1$, such that $(S_1,{\frac{1}{c}}T_1)$ is  $m_1$-symmetric and $(S_2,cT_2)$ is  $m_2$-symmetric.
\end{pro}
Looking upon $\delta^m_{S_1\otimes S_2,T_1\otimes T_2}(I\otimes I)$ as the sum
\begin{eqnarray*}
\delta^m_{S_1\otimes S_2,T_1\otimes T_2}(I\otimes I)&=&  \left(\sum_{j=0}^{m} (-1)^j\left(\begin{array}{clcr}m\\j\end{array}\right)L^{m-j}_{S_1\otimes S_2}R^j_{T_1\otimes T_2}\right)(I\otimes I)\\
&=&  \sum_{j=0}^{m}(-1)^j \left(\begin{array}{clcr}m\\j\end{array}\right) S^{m-j}_1T^j_1\otimes S^{m-j}_2T^j_2\\
&=&  \sum_{j=0}^{m} (-1)^j\left(\begin{array}{clcr}m\\j\end{array}\right)\left(\lfloor S_1\rfloor^{m_i-j} \times \lfloor T_1\rfloor^j\right) \otimes \left(\lfloor S_2\rfloor^{m-j} \times \lfloor T_2\rfloor^j\right),
\end{eqnarray*}
Proposition \ref{pro01} says the following.
\begin{pro}\label{pro11} If $S_1, S_2\in\b$ are invertible operators such that
$$
 \sum_{j=0}^{m}(-1)^j \left(\begin{array}{clcr}m\\j\end{array}\right)\left(\lfloor S_1\rfloor^{m-j} \times \lfloor T_1\rfloor^j\right) \otimes \left(\lfloor S_2\rfloor^{m-j} \times \lfloor T_2\rfloor^j\right)=0,
$$
and
$$
 \sum_{j=0}^{m-1}(-1)^j \left(\begin{array}{clcr}m-1\\j\end{array}\right)\left(\lfloor S_1\rfloor^{m-1-j} \times \lfloor T_1\rfloor^j\right) \otimes \left(\lfloor S_2\rfloor^{m-1-j} \times \lfloor T_2\rfloor^j\right)\neq 0,
$$
then there exists a non-zero scalar $c$ and positive integers $m_i$, $m=m_1+m_2-1$, such that 
$$
 \sum_{j=0}^{m_1}(-1)^j \left(\begin{array}{clcr}m_1\\j\end{array}\right)\left(\lfloor S_1\rfloor^{m_1-j} \times \lfloor {\frac{1}{c}}T_1\rfloor^j\right)=0= \sum_{j=0}^{m_2} (-1)^j\left(\begin{array}{clcr}m_2\\j\end{array}\right)\left(\lfloor S_2\rfloor^{m-j} \times \lfloor cT_2\rfloor^j\right)
$$
and 
\begin{eqnarray*}
& & \sum_{j=0}^{m_1-1}(-1)^j \left(\begin{array}{clcr}m_1-1\\j\end{array}\right)\left(\lfloor S_1\rfloor^{m_1-1-j} \times \lfloor {\frac{1}{c}}T_1\rfloor^j\right)\neq 0\\
&\neq& \sum_{j=0}^{m_2-1}(-1)^j \left(\begin{array}{clcr}m_2-1\\j\end{array}\right) \left(\lfloor S_2\rfloor^{m_2-1-j} \times \lfloor cT_2\rfloor^j\right).
\end{eqnarray*}
\end{pro}
Corresponding to Theorem \ref{thm12}, we have the following result for tensor products of commuting $d$-tuples satisfying a strict $m$-symmetric property.

\begin{thm}\label{thm22} If $\S_i, \T_i\in\b^d$, $1\leq i\leq 2$, are commuting $d$-tuples, $\S_1$ and  $\S_2$ are invertible, and $(\S_1\otimes\S_2,\T_1\otimes\T_2)$ is strict $m$ symmetric, 
then there exists a non-zero scalar $c$ and positive integers $m_i$, $m=m_1+m_2-1$, such that $(\S_1,{\frac{1}{c}}\T_1)$ is strict $m_1$-symmetric and $(\S_2,c\T_2)$ is strict $m_2$-symmetric.
\end{thm}
\begin{demo} If $(\S_1\otimes\S_2,\T_1\otimes\T_2)$ is strict $m$-symmetric, then
\begin{eqnarray*}
& & \delta^m_{\S_1\otimes\S_2,\T_1\otimes\T_2}(I\otimes I)\\
&=& \left(\sum_{j=0}^{m}(-1)^j \left(\begin{array}{clcr}m\\j\end{array}\right)\L^{m-j}_{\S_1\otimes\S_2} \times \R^j_{\T_1\otimes\T_2}\right)(I\otimes I)\\
&=& \sum_{j=0}^{m}(-1)^j \left(\begin{array}{clcr}m\\j\end{array}\right)\left(\sum_{i,k=1}^d S_{1i} \otimes S_{2k}\right)^{m-j}\left( \sum_{i,k=1}^d T_{1i} \otimes T_{2k}\right)^j\\
&=& \sum_{j=0}^{m}(-1)^j \left(\begin{array}{clcr}m\\j\end{array}\right)\left(\lfloor\d S_{1i}\rfloor^{m-j}  \lfloor  \d T_{1i}\rfloor^j\right)  \otimes \left( \lfloor \d S_{2i}\rfloor^{m-j}   \lfloor \d T_{2i}\rfloor^j\right)\\
&=& 0,\\
&{\rm and}&\\
& & \delta^{m-1}_{\S_1\otimes\S_2,\T_1\otimes\T_2}(\I\otimes\I)\\
&=&  \sum_{j=0}^{m-1}(-1)^j \left(\begin{array}{clcr}m-1\\j\end{array}\right)\left(\lfloor\d S_{1i}\rfloor^{m-1-j}  \lfloor  \d T_{1i}\rfloor^j\right)  \otimes \left( \lfloor \d S_{2i}\rfloor^{m-1-j}  \lfloor \d T_{2i}\rfloor^j\right)\\
&\neq& 0.
\end{eqnarray*}
The operators $\S_1$ and $\S_2$ being invertible, $\d S_{1i}$ and $\d S_{2i}$ are invertible, Proposition \ref{pro11} applies and we conclude the existence of a non-zero scalar $c$ and positive integers $m_i$, $m=m_1+m_2-1$, such that

\begin{eqnarray*}
& &  \sum_{j=0}^{m_1}(-1)^j \left(\begin{array}{clcr}m_1\\j\end{array}\right)\left(\lfloor\d S_{1i}\rfloor^{m_1-j}  \lfloor \d {\frac{1}{c}}T_{1i}\rfloor^j\right)=\delta^{m_1}_{\S_1,{\frac{1}{c}}\T_1}(\I)=0,\\
& &  \sum_{j=0}^{m_2}(-1)^j \left(\begin{array}{clcr}m_2\\j\end{array}\right)\left(\lfloor\d S_{2i}\rfloor^{m-j}  \lfloor  \d cT_{2i}\rfloor^j\right)=\delta^{m_2}_{\S_2,c\T_2}(\I)=0
\end{eqnarray*}
and
$$
\delta^{m_1-1}_{\S_1,{\frac{1}{c}}\T_1}(\I)\neq 0,  \  \delta^{m_2-1}_{\S_2,c\T_2}(\I)\neq 0.
$$
This completes the proof.
\end{demo}
\begin{rema}\label{rema22} {\em Paul and Gu \cite[Theorem 5.2]{PG} state that ``if the operators $S_i$ are left invertible and the operators $T_i$ are right invertible, $1\leq i\leq 2$, then $(S_1\otimes S_2, T_1\otimes T_2)$ is $m$-symmetric if and only if there exist a non-zero scalar $c$ and positive integers $m_i$, $m=m_1+m_2-1$, such that $(S_1,{\frac{1}{c}} T_1)$ is strict $m_1$-symmetric and $(S_2,c T_2)$ is strict $m_2$-symmetric". The hypothesis $S_i$ are left invertible and $T_i$ are invertible is a bit of an overkill, as we show below. As seen in the proof of Theorem \ref{thm22}, the invertibility of $S_1$ and $S_2$  - a fact guraranteed by the left invertibility of $S_i$ and the right invertibility of $T_i$ - is sufficient. If $S_i$, $1\leq i\leq 2$, is left invertible, then there exist operators $E_i$ such that $(E_1\otimes E_2)(S_1\otimes S_2)=(I\otimes I)$ ($=\I$) and 
\begin{eqnarray*}
& & 0=\delta^m_{S_1\otimes S_2,T_1\otimes T_2}(\I)\\
&=& \sum_{j=0}^{m}(-1)^j \left(\begin{array}{clcr}m\\j\end{array}\right)(S_1\otimes S_2)^{m-j}(T_1\otimes T_2)^j\\
& = & \sum_{j=0}^{m}(-1)^j \left(\begin{array}{clcr}m\\j\end{array}\right)(E_1\otimes E_2)^j(T_1\otimes T_2)^j\\
&=& \triangle^m_{E_1\otimes E_2,T_1\otimes T_2}(\I).
\end{eqnarray*}
For conveneience, set $E_1\otimes E_2=A$ and $T_1\otimes T_2=B$. Then $\triangle^m_{A,B}(\I)=0$. It is easily seen, use induction, that $(a-1)^m=a^m-\sum_{j=0}^{m-1} \left(\begin{array}{clcr}m\\j\end{array}\right)(a-1)^j$; hence
$$
(L_AR_B-\I)^m=(L_AR_B)^m - \sum_{j=0}^{m-1} \left(\begin{array}{clcr}m\\j\end{array}\right)(L_AR_B-\I)^j
$$
and upon letting $(L_AR_B-\I)=\nabla_{A,B}$ that
\begin{eqnarray*}
& & \nabla^m_{A,B}(\I)=0\Longleftrightarrow (L_AR_B)^m(\I) - \sum_{j=0}^{m-1} \left(\begin{array}{clcr}m\\j\end{array}\right)\nabla_{A,B}^j(\I)=0\\
&\Longrightarrow&  (L_AR_B)^{m+1}(\I) = \sum_{j=0}^{m-1} \left(\begin{array}{clcr}m\\j\end{array}\right)(L_AR_B)\nabla_{A,B}^j(\I)\\
&=& \sum_{j=0}^{m-1} \left(\begin{array}{clcr}m\\j\end{array}\right)\nabla_{A,B}^{j+1}(\I) +\sum_{j=0}^{m-1} \left(\begin{array}{clcr}m\\j\end{array}\right)\nabla_{A,B}^j(\I)\\
&=& \left(\begin{array}{clcr}m\\m-1\end{array}\right)\nabla^m_{A,B}(\I) +\sum_{j=0}^{m-1} \left(\begin{array}{clcr}m+1\\j\end{array}\right)\nabla_{A,B}^j(\I)\\
&=& \sum_{j=0}^{m-1} \left(\begin{array}{clcr}m+1\\j\end{array}\right)\nabla_{A,B}^j(\I).
\end{eqnarray*}
An induction argument now proves
\begin{eqnarray*} 
& & (L_A R_B)^n(\I) =\sum_{j=0}^{m-1} \left(\begin{array}{clcr}n\\j\end{array}\right)\nabla_{A,B}^j(\I)\\
&=&  \left(\begin{array}{clcr}n\\m-1\end{array}\right)\nabla^{m-1}_{A,B}(\I) + \sum_{j=0}^{m-2} \left(\begin{array}{clcr}n\\j\end{array}\right)\nabla_{A,B}^j(\I)
\end{eqnarray*}
for all integers $n\geq m$. Observe that $\triangle^m_{A,B}(\I)=0$ implies $A$ is right invertible and $B$ is left invertible; since already $B=T_1\otimes T_2$ is right invertible, $B$ is invertible, and then
\begin{eqnarray*}
\delta^m_{A,B}(\I)=0 &\Longleftrightarrow& \sum_{j=0}^{m}(-1)^j \left(\begin{array}{clcr}m\\j\end{array}\right)L_A^{m-j}R_B^{j}(\I)=0\\
&\Longleftrightarrow& \sum_{j=0}^{m}(-1)^j \left(\begin{array}{clcr}m\\j\end{array}\right)L_A^{m-j}R_B^{-m+j}(\I)=0\\
&\Longleftrightarrow& \triangle^m_{A,B^{-1}}(\I)=0\Longrightarrow A \ {\rm is \ right \ invertible}\Longrightarrow A \ {\rm is \ invertible}.
\end{eqnarray*}
The invertibility of  $A$ and  $B$ implies that of $L_AR_B$. We have
\begin{eqnarray*}
{\frac{1}{  \left(\begin{array}{clcr}n\\m-1\end{array}\right)}}\left(\I - \sum_{j=0}^{m-2} \left(\begin{array}{clcr}n\\j\end{array}\right)(L_AR_B)^{-n}\nabla_{A,B}^j(\I)\right)=\nabla^{m-1}_{A,B}(\I).
\end{eqnarray*}
Since $ \left(\begin{array}{clcr}n\\m-1\end{array}\right)$ is of the order of $n^{m-1}$ and $ \left(\begin{array}{clcr}n\\m-2\end{array}\right)$ is of the order of $n^{m-2}$, letting $n\longrightarrow\infty$ this implies 
$$\nabla^{m-1}_{A,B}(\I)=0\Longleftrightarrow \triangle^{m-1}_{A,B}(\I)=0.$$
Repeating the argument, we eventually have
$$
\triangle_{A,B}(\I)=0\Longleftrightarrow (S^{-1}_1\otimes S^{-1}_2)(T_1\otimes T_2)=\I\Longleftrightarrow S_1\otimes S_2=T_1\otimes T_2;
$$
hence  there exists a scalar $c$ such that $S_1=cT_1$ and $S_2={\frac{1}{c}}T_2$. In particular, if $S_i, T_i$ are Hilbert space operators such that $S_i=T_i^*$, then $T_1\otimes T_2$ is self-adjoint.
}\end{rema}


\vskip10pt \noindent\normalsize\rm B.P. Duggal,{University of Ni\v s,
Faculty of Sciences and Mathematics,
P.O. Box 224, 18000 Ni\v s, Serbia}.

\noindent\normalsize \tt e-mail:  bpduggal@yahoo.co.uk

\end{document}